\documentclass{article}

\usepackage{latexsym}

\newtheorem{llemma}{Lemma}[section]
\newtheorem{prop}[llemma]{Proposition}
\newtheorem{exmp}[llemma]{Example}

\newtheorem{ccorollary}[llemma]{Corollary}
\newtheorem{defn}[llemma]{Definition}

\newtheorem{key}[llemma]{Keyword}

\begin{document}

%*********************************************
% TITLE

\title{Relations on FP-Soft Sets Applied to Decision Making Problems}

\author{Irfan Deli\\
Department of Mathematics, Faculty of Arts and Sciences\\
7 Aral{\i}k University, 79000 Kilis, Turkey \\
irfandeli@kilis.edu.tr
\\~\\
Naim \c{C}a\u{g}man \\
Department of Mathematics, Faculty of Arts and Sciences\\
Gaziosmanpa\c{s}a University, 60250 Tokat, Turkey\\
naim.cagman@gop.edu.tr }

\maketitle

%**********************************************
% ABSTRACT

\begin{abstract}
In this work, we first define relations on the fuzzy parametrized
soft sets and study their properties. We also give a decision making
method based on these relations. In approximate reasoning, relations
on the fuzzy parametrized soft sets have shown to be of a primordial
importance. Finally, the method is successfully applied to a
problems that contain uncertainties.
\begin{key}
Soft sets, fuzzy sets, FP-soft sets, relations on FP-soft sets,
decision making.
\end{key}
\end{abstract}

%**********************************************
% SECTION 1: INTRODUCTION

\section{Introduction}

In 1999, the concept of soft sets was introduced by Molodtsov
\cite{mol-99-teo} to deal with problems that contain uncertainties.
After Molodtsov, the operations of soft sets are given in
\cite{cag-10a-teo, maj-03-teo, op} and studied their properties.
Since then, based on these operations, soft set theory has developed
in many directions and applied to wide variety of fields. For
instance; on the theory of soft sets
\cite{ali-09-teo,cag-10a-teo,cag-10b-teo,cag-11f,kar-12,maj-03-teo,maju-10-teo,op},
on the soft decision making
\cite{fee,fen-10-dec,fen-12-dec,kov-07-dec,maj-02-dec,roy-07-dec},
on the fuzzy soft sets \cite{cag-12b,cag-10-pfss,cag-10-fss} and
soft rough sets \cite{fee} are some of the selected works. Some
authors have also studied the algebraic properties of soft sets,
such as
\cite{aca-10-alg,akt-07,cag-12a,jun-10-alg,par-08-alg,op2,op3}.

The FP-soft sets, firstly studied by \c{C}a\u{g}man \emph{et al.}
\cite{cag-11e}, is a fuzzy parameterized soft sets. Then, FP-soft
sets theory and its applications studied in detail, for example
\cite{deli-12a,deli-12b,deli-12}. In this paper, after given most of
the fundamental definitions of the operations of fuzzy sets, soft
sets and FP-soft sets in next section, we define relations on
FP-soft sets and we also give their properties in Section 3. In
Section 4, we define symmetric, transitive and reflexive relations
on the FP-soft sets. In Section 5, we construct a decision making
method based on the FP-soft sets. We also give an application which
shows that this methods successfully works. In the final section,
some concluding comments are presented.
\section{Preliminary}\label{ss}
%-----------------------------------------------------------
In this section, we give the basic definitions and results of soft
set theory \cite{mol-99-teo} and fuzzy set theory \cite{zad-65} that
are useful for subsequent discussions.
%-----------------------------------------------------------------
\begin{defn}\cite{zad-65}
Let $U$ be the universe. A fuzzy set $X$ over $U$ is a set defined
by a membership function $\mu_X$ representing a mapping
$$
\mu_X:U\rightarrow [0,1].
$$
The value $\mu_X(x)$ for the fuzzy set
$X$ is called the membership value or the grade of membership of
$x\in U$. The membership value represents the degree of $x$
belonging to the fuzzy set $X$. Then a fuzzy set $X$ on $U$ can be
represented as follows, $$X=\{(\mu_X(x)/x):x\in U,\mu_X(x)\in
[0,1]\}.$$
\end{defn}
Note that the set of all fuzzy sets on $U$ will be denoted by
$F(U)$.
%--------------------------------------------------
\begin{defn}\label{add}\cite{dub-80}
$t$-norms are associative, monotonic and commutative two valued
functions $t$ that map from $[0,1]\times [0,1]$ into $[0,1]$. These
properties are formulated with the following conditions:
%------------------------
\begin{enumerate}
\item $t(0,0)=0$ and $t(\mu_{X_1}(x),1)=t(1,\mu_{X_1}(x))=\mu_{X_1}(x),\,\, x\in E$
\item If  $\mu_{X_1}(x)\leq \mu_{X_3}(x)$ {and} $\mu_{X_2}(x)\leq \mu_{X_4}(x)$, then\\
     $t(\mu_{X_1}(x),\mu_{X_2}(x))\leq t(\mu_{X_3}x),\mu_{X_4}(x))$
\item $t(\mu_{X_1}(x),\mu_{X_2}(x))= t(\mu_{X_2}(x),\mu_{X_1}(x))$
\item $t(\mu_{X_1}(x),t(\mu_{X_2}(x),\mu_{X_3}(x)))=t(t(\mu_{X_1}(x),\mu_{X_2})(x),\mu_{X_3}(x))$
\end{enumerate}
%------------------------
\end{defn}
%----------------------------------------------------
\begin{defn}\label{add}\cite{dub-80}
$t$-conorms or $s$-norm are associative, monotonic and commutative
two placed functions $s$ which map from $[0,1]\times [0,1]$ into
$[0,1]$. These properties are formulated with the following
conditions:
\begin{enumerate}
\item  $s(1,1)=1 \,\,and \,\, s(\mu_{X_1}(x),0)=s(0,\mu_{X_1}(x))=\mu_{X_1}(x), x\in E$
\item  {if}
$\mu_{X_1}(x)\leq \mu_{X_3}(x)$ {and}  $\mu_{X_2}(x)\leq
\mu_{X_4}(x)$, then \\
$s(\mu_{X_1}(x),\mu_{X_2}(x))\leq s(\mu_{X_3}(x),\mu_{X_4}(x))$
\item   $s(\mu_{X_1}(x),\mu_{X_2}(x))= s(\mu_{X_2}(x),\mu_{X_1}(x))$
\item  $s(\mu_{X_1}(x),s(\mu_{X_2}(x),\mu_{X_3}(x)))=s(s(\mu_{X_1}(x),\mu_{X_2})(x),\mu_{X_3}(x))$
\end{enumerate}
\end{defn}
%----------------------------------------------------
$t$-norm and $t$-conorm are related in a sense of lojical duality.
Typical dual pairs of non parametrized $t$-norm and $t$-conorm are
complied below:

\begin{enumerate}
\item  Drastic product:
$$
t_w(\mu_{X_1}(x),\mu_{X_2}(x)) = \left \{\begin{array}{ll}
min\{\mu_{X_1}(x),\mu_{X_2}(x)\}, & max\{\mu_{X_1}(x)  \mu_{X_2}(x)\}=1\\
0,                                & otherwise \\
\end{array}\right.
$$
\item  Drastic sum:
$$
s_w(\mu_{X_1}(x),\mu_{X_2}(x)) = \left \{\begin{array}{ll}
max\{\mu_{X_1}(x),\mu_{X_2}(x)\}, &  min\{\mu_{X_1}(x)  \mu_{X_2}(x)\}=0\\
1, & otherwise\\
\end{array}\right.
$$
\item  Bounded product:
$$
t_1(\mu_{X_1}(x),\mu_{X_2}(x))=max\{0, \mu_{X_1}(x)+\mu_{X_2}(x)-1\}
$$
\item  Bounded sum:
$$
s_1(\mu_{X_1}(x),\mu_{X_2}(x))=min\{1, \mu_{X_1}(x)+\mu_{X_2}(x)\}
$$
\item Einstein product:
 $$
 t_{1.5}(\mu_{X_1}(x),\mu_{X_2}(x))=\frac{\mu_{X_1}(x).\mu_{X_2}(x)}{2-[\mu_{X_1}(x)+\mu_{X_2}(x)-\mu_{X_1}(x).\mu_{X_2}(x)]}
 $$
\item  Einstein sum:
$$
s_{1.5}(\mu_{X_1}(x),\mu_{X_2}(x))=\frac{\mu_{X_1}(x)+\mu_{X_2}(x)}{1+\mu_{X_1}(x).\mu_{X_2}(x)}
$$
\item  Algebraic product:
$$
t_{2}(\mu_{X_1}(x),\mu_{X_2}(x))=\mu_{X_1}(x).\mu_{X_2}(x)
$$
\item  Algebraic sum:
$$
s_{2}(\mu_{X_1}(x),\mu_{X_2}(x))=\mu_{X_1}(x)+\mu_{X_2}(x)-\mu_{X_1}(x).\mu_{X_2}(x)
$$
\item  Hamacher product: $$
t_{2.5}(\mu_{X_1}(x),\mu_{X_2}(x))=\frac{\mu_{X_1}(x).\mu_{X_2}(x)}{\mu_{X_1}(x)+\mu_{X_2}(x)-\mu_{X_1}(x).\mu_{X_2}(x)}
$$
\item  Hamacher sum:
$$
s_{2.5}(\mu_{X_1}(x),\mu_{X_2}(x))=\frac{\mu_{X_1}(x)+\mu_{X_2}(x)-2.\mu_{X_1}(x).\mu_{X_2}(x)}{1-\mu_{X_1}(x).\mu_{X_2}(x)}
$$
\item  Minumum:
$$
t_3(\mu_{X_1}(x),\mu_{X_2}(x))=min\{\mu_{X_1}(x),\mu_{X_2}(x)\}
$$
\item  Maximum:
$$
s_3(\mu_{X_1}(x),\mu_{X_2}(x))=max\{\mu_{X_1}(x),\mu_{X_2}(x)\}
$$
\end{enumerate}
%-------------------------------------------------------------
\begin{defn}\cite{mol-99-teo}.
Let $U$ be an initial universe set and let E be a set of parameters.
Then, a pair $(F,E)$ is called a soft set over $U$ if and only if
$F$ is a mapping or $E$ into the set of aft subsets of the set $U$.

In other words, the soft set is a parametrized family of subsets of
the set $U$. Every set $F(\varepsilon)$, $\varepsilon\in E$, from
this family may be considered as the set of $\varepsilon$-elements
of the soft set $(F, E)$, or as the set of $\varepsilon$-approximate
elements of the soft set.

It is worth noting that the sets $F(\varepsilon)$ may be arbitrary.
Some of them may be empty, some may have nonempty intersection.
\end{defn}
%-------------------------------------------------------------
In this definition, $E$ is a set of parameters that are describe the
elements of the universe $U$. To apply the soft set in decision
making subset $A, B, C,...$ of the parameters set $E$ are needed.
Therefore, \d{C}a\u{g}man and Engino\u{g}lu \cite{cag-10a-teo}
modified the definition of soft set as follows.
%-------------------------------------------------------------
\begin{defn}\cite{cag-10a-teo}
Let $U$ be a universe, $E$ be a set of parameters that are describe
the elements of $U$, and $A\subseteq E$. Then, a soft set $F_A$ over
$U$ is a set defined by a set valued function $f_A$ representing a
mapping
%---------
\begin{equation}\label{soft-set}
f_A: E\to P(U) \textrm{ such that}\, f_A(x)=\emptyset \textrm{ if }
x\in E-A
\end{equation}
%---------
where $f_A$ is called approximate function of the soft set $F_A$. In
other words, the soft set is a parametrized family of subsets of the
set $U$, and therefore it can be written a set of ordered pairs
$$
F_A= \{(x, f_A(x)): x\in E, f_A(x)=\emptyset \textrm{ if } x\in
E-A\}
$$
\end{defn}
The subscript $A$ in the $f_A$ indicates that $f_A$ is the
approximate function of $F_A$. The value $f_A(x)$ is a set called
\emph{$x$-element} of the soft set for every $x\in E$.
%
%Note that, from now on, the set of all soft sets over $U$ will be
%denoted by $S(U)$.
%-----------------------------------------------------------
\begin{defn} \cite{cag-11e}
Let $F_X$ be a soft set over $U$ with its approximate function $f_X$
and $X$ be a fuzzy set over $E$ with its membership function
$\mu_X$. Then, a $FP-$soft sets $\Gamma_X$, is a fuzzy parameterized
soft set over $U$, is defined by the set of ordered pairs
$$
\label{FP-soft set} \Gamma_X= \{(\mu_X(x)/x, f_X(x)): x\in E\}
$$
where $f_X: E\to P(U)$ such that $f_X(x)=\emptyset$ if $\mu_X(x)=0$
is called approximate function and $\mu_X:E\to [0,1]$ is called
membership function of $FP-$soft set $\Gamma_X$. The value
$\mu_X(x)$ is the degree of importance of the parameter $x$ and
depends on the decision-maker's requirements.
\end{defn}
%-----------------------------------------------------------
Note that the sets of all $FP-$soft sets over $U$ will be denoted by
$FPS(U)$.

%------------------------------------------------------------
%\begin{defn}\label{sfo-or} \cite{deli-12}
%Let $\Gamma_X\in FPS(U)$, $f_{X}$ be an approximation function and
%$\mu_{X}$ be a membership function of $\Gamma_X$. Then, an $FP-$soft
%fuzzification operator, denoted $s_{\Gamma_X}$, is defined by
%$$
%s_{\Gamma_X}: F(E)\times U\to F(U),
%$$
%$$
%s_{\Gamma_X}(\mu_X(x_i),u_j)=\{\mu_X(u_j)/u_j: u_j\in U,
%\mu_X(u_j)\in [0,1],i=1,2,...,n\}
%$$
%where
%$$
%\mu_X(u_j)=\frac{1}{|U|}\sum_{i}\mu_X(x_i)\chi(u_j)
%$$
%and where
%$$
%\chi(u_j)=\left\{
%\begin{array}{ll}
%1, & u_j\in f_X(x_i)\\
%0, & u_j \notin f_X(x_i)
%\end{array}\right.
%$$
%
%\end{defn}
%------------------------------------------------------------
\section{Relations on the FP-Soft Sets}
In this section, after given the cartesian products of two FP-soft
sets, we define a relations on FP-soft sets and study their desired
properties.
\begin{defn}
Let $\Gamma_X,\Gamma_Y\in FPS(U)$. Then, a cartesian product of
$\Gamma_X$ and $\Gamma_Y$, denoted by $\Gamma_X\widehat{{\times}}
\Gamma_Y$, is defined as
$$
\Gamma_X\widehat{{\times}}
\Gamma_Y=\Big\{\big(\mu_{X\widehat{{\times}}
Y}(x,y)/(x,y),f_{X\widehat{{\times}} Y}(x,y)):(x,y)\in E{\times}
E\big)\Big\}$$ where
$$
f_{X\widehat{{\times}} Y}(x,y)=f_X(x)\cap f_Y(y)
$$
and
$$
\mu_{X\widehat{{\times}} Y}(x,y)=min\{\mu_X(x), \mu_Y(y)\}
$$
\end{defn}
Here $\mu_{X\widehat{{\times}} Y}(x,y)$ is a t-norm.

%----------------------------------------------------------------------------------
\begin{exmp}\label{1}
Let $U=\{u_1, u_2, u_3, u_4, u_5, u_6, u_7, u_8, u_9, u_{10},
u_{11}, u_{12}, u_{13}, u_{14}, u_{15}\}$, $E=\{x_1, x_2, x_3, x_4,
x_5, x_6, x_7, x_8\}$, and  $X= \{0.5/x_1, 0.7/x_2, 0.3/x_3,
0.9/x_4, 0.6/x_5\}$ and $Y=\{0.9/x_3, 0.1/x_6, 0.7/x_7, 0.3/x_8\}$
be two fuzzy subsets of $E$. Suppose that
$$
\begin{array}{rl}
\Gamma_X= & \bigg\{(0.5/x_1,\{u_1, u_3, u_4,u_6, u_7,u_8,
 u_{11},u_{12},u_{13},u_{15}\}),(0.7/x_2,\{u_3, u_7, u_8,u_{14},\\&u_{15}\}),
(0.3/x_3,\{u_1,u_2,u_4,u_5, u_6,u_9,u_{10},u_{12},u_{13}\}),
(0.9/x_4,\{u_2,u_4,u_6,u_8,\\&u_{12},u_{13}\}),
(0.6/x_5,\{u_3,u_4,u_6,u_7,u_9,u_{13},u_{15}\})\bigg\} \\
\Gamma_Y= & \bigg\{( 0.9/x_3,\{u_1,u_5, u_6,u_9, u_{10},
u_{13}\}),(0.1/x_6,\{u_3,u_5, u_7, u_8,u_9,u_{11},u_{15}\}),\\&
(0.7/x_7,\{u_2,u_5,u_9, u_{10},u_{11}, u_{14}\}),(0.3/x_8,\{u_2,
u_5,u_8,u_{10},u_{12},u_{14}\})\bigg\} \\
\end{array}
$$
Then, the cartesian product of $\Gamma_X$ and $\Gamma_Y$ is obtained
as follows
$$
\begin{array}{rl}

\Gamma_X\widehat{\times} \Gamma_Y=& \bigg\{
(0.5/(x_1,x_3),\{u_1,u_6, u_{13}\}),(0.1/(x_1,x_6),\{u_3, u_7,
u_8,u_{11},u_{15}\}),\\&
(0.5/(x_1,x_7),\{u_{11}\}),(0.3/(x_1,x_8),\{u_8,u_{12}\})
,(0.7/(x_2,x_3),\emptyset),\\&
 (0.1/(x_2,x_6),\{u_3,u_7,u_8\}),
(0.7(x_2,x_7),\{u_{14}\}), (0.3/(x_2,x_8),\\&\{u_8,u_{14}\}),
(0.3/(x_3,x_3),\{u_1,u_5, u_6,u_9,u_{10}, u_{13}\}),
 (0.1/(x_3,x_6),\\&\{u_5,u_9\}),
 ( 0.3/(x_3,x_7),\{u_2,u_5,u_9,u_{10}\}),
 (0.3/(x_3,x_8),\\&\{u_2,u_5,u_{10},u_{12}\}),
(0.9/(x_4,x_3),\{u_6\}), (0.1/(x_4,x_6),\emptyset),\\&
(0.7/(x_4,x_7),\{u_2,u_6\}),
 (0.3/(x_4,x_8),\{u_2,u_8,u_{12}\}),\\&
(0.6/(x_5,x_3),\{u_6,u_9,u_{13}\}),(0.1/(x_5,x_6),\{u_3,u_7,u_9,u_{11},u_{15}\}),
\\&(0.6/(x_5,x_7),\{u_2\}),(0.3/(x_5,x_8),\emptyset)\bigg\} \\
\end{array}
$$
\end{exmp}
%---------------------------------------------------------
\begin{defn}
Let $\Gamma_X,\Gamma_Y\in FPS(U)$. Then, an FP-soft relation from
$\Gamma_X$ to $\Gamma_Y$, denoted by $R_F$,  is an FP-soft subset of
$\Gamma_X\widehat{\times}\Gamma_Y$. Any FP-soft subset of
$\Gamma_X{\times}\Gamma_Y$ is called a FP-relation on $\Gamma_X$.
\end{defn}
Note that if $\alpha=(\mu_X(x),f_X(x))\in \Gamma_X$ and
$\beta=(\mu_Y(y),f_Y(y))\in \Gamma_Y$, then
$$
\alpha R_F\beta\Leftrightarrow \big(\mu_{X\widehat{\times}
Y}(x,y)/(x,y),f_{X\widehat{\times} Y}(x,y))\in R_F
$$

%----------------------------------------------------------------------------------
\begin{exmp}\label{2}Let us consider the Example \ref{1}. Then, we define
an FP-soft relation $R_F$, from $\Gamma_Y$ to $\Gamma_X$, as follows
$$
\alpha{R_{F}}\beta\Leftrightarrow \mu_{X\widehat{\times}
Y}(x_i,x_j)/(x_i,x_j))\geq0.3\quad(1\leq i,j \leq 3)
$$
Then
$$
\begin{array}{rl}
 {R_F}=& \bigg\{
(0.5/(x_1,x_3),\{u_1,u_6, u_{13}\}),
((0.5/(x_1,x_7),\{u_{11}\}),(0.3/(x_1,x_8),\{u_8,\\&u_{12}\}),
(0.7(x_2,x_7),\{u_{14}\}), (0.3/(x_2,x_8),\{u_8,u_{14}\}),
(0.3/(x_3,x_3),\{u_1,\\&u_5, u_6,u_9,u_{10}, u_{13}\}),
 ( 0.3/(x_3,x_7),\{u_2,u_5,u_9,u_{10}\}),
 (0.3/(x_3,x_8),\{u_2,\\&u_5,u_{10},u_{12}\}),
(0.9/(x_4,x_3), \{u_6\}),(0.7/(x_4,x_7),\{u_2,u_6\}),
 (0.3/(x_4,x_8), \\&\{u_2,u_8,u_{12}\}),
(0.6/(x_5,x_3),\{u_6,u_9,u_{13}\}),
(0.6/(x_5,x_7),\{u_2\})\bigg\} \\
\end{array}
$$
\end{exmp}
%----------------------------------------------------------------------------------
%\begin{exmp}\label{3} Let us consider the Example \ref{2}. Then, we define
%two FP-soft relation $R_{F_1}$ and $R_{F_1}$ from $\Gamma_X$ to
%$\Gamma_Y$ as follows
%$$
%\alpha{R_{F_1}}\beta\Leftrightarrow \mu_{X\widehat{\times}
%Y}(x_i,x_j)/(x_i,x_j))\geq0.5\quad(i\in\{1,2,3,4,5\} \,\,and \,\,
%j\in\{3,6,7,8\} )
%$$
%Then
%$$
%\begin{array}{rl}
% {R_{F_1}}=& \bigg\{(0.5/(x_2,x_3),\{u_1,u_6,u_{10}, u_{13}\}),
%(0.5/(x_1,x_7)=\{u_{11}\}) ,
%\\&(0.7(x_2,x_7)=\{u_{14}\}),(0.9/(x_4,x_3)= \{u_6\}),
%\\&(0.7/(x_4,x_7)=\{u_2,u_6\}),(0.6/(x_5,x_3)=\{u_6,u_9, u_{13}\}),
%\\&(0.6/(x_5,x_7)=\{u_2\})\bigg\} \\
%\end{array}
%$$
%and,
%$$
%\alpha{R_{F_2}}\beta\Leftrightarrow \mu_{X\widehat{\times}ð
%Y}(x_i,x_j)/(x_i,x_j))\geq0.7\quad(i\in\{1,2,3,4,5\} \,\,and \,\,
%j\in\{3,6,7,8\} )
%$$
%Then
%$$
%\begin{array}{rl}
% {R_{F_2}}=& \bigg\{(0.7/(x_2,x_3)=\emptyset),
%(0.7(x_2,x_7)=\{u_{14}\}),((0.9/(x_4,x_3)= \{u_6\}),\\&
%(0.7/(x_4,x_7)=\{u_2,u_6\})\bigg\} \\
%\end{array}
%$$
%\end{exmp}
%-----------------------------------------------------------------------------
\begin{defn}
Let $\Gamma_X,\Gamma_Y\in FPS(U)$ and  $ R_F$ be an FP-soft relation
from $\Gamma_X$ to $\Gamma_Y$. Then  domain and range of $ R_F$
respectively is defined as
\begin{eqnarray*}
D( R_F) & = & \{\alpha\in F_A:\alpha R_F\beta\}\\
R( R_F) & = & \{\beta\in F_B:\alpha R_F\beta\}.
\end{eqnarray*}
\end{defn}
%----------------------------------------------------------------------------------

\begin{exmp}Let us consider the Example \ref{2}.
$$
\begin{array}{rl}
D( {R_{F}})= & \bigg\{(0.5/x_1,\{u_1, u_3, u_4,u_6, u_7,u_8,
 u_{11},u_{12},u_{13},u_{15}\}),(0.7/x_2,\{u_3, u_7, \\&u_8,u_{14},u_{15}\}),
(0.3/x_3,\{u_1,u_2,u_4,u_5, u_6,u_9,u_{10},u_{12},u_{13}\}),
(0.9/x_4,\\&\{u_2,u_4,u_6,u_8,u_{12},u_{13}\}),
(0.6/x_5,\{u_3,u_4,u_6,u_7,u_9,u_{13},u_{15}\})\bigg\} \\
\end{array}
$$
$$
\begin{array}{rl}
 R({R_{F}})= & \bigg\{( 0.9/x_3,\{u_1,u_5, u_6,u_9, u_{10},
u_{13}\}), (0.7/x_7,\{u_2,u_5,\\&u_9, u_{10},u_{11},
u_{14}\}),(0.3/x_8,\{u_2,
u_5,u_8,u_{10},u_{12},u_{14}\})\bigg\} \\
\end{array}
$$
\end{exmp}
%---------------------------------------------------------------------------------
\begin{defn}
Let $R_F$ be an FP-soft relation from $\Gamma_X$ to $\Gamma_Y$. Then
$R_F^{-1}$ is  from $\Gamma_Y$ to $\Gamma_X$ is defined as
$$\alpha R_F^{-1}\beta= \beta R_F\alpha$$
\end{defn}
\begin{exmp}\label{4}Let us consider the Example \ref{2}.
Then, ${R_{F}}^{-1}$ is  from $\Gamma_Y$ to $\Gamma_X$ is obtained
by
$$
\begin{array}{rl}
{R_{F}}^{-1}=& \bigg\{ (0.5/(x_3,x_1),\{u_1,u_6, u_{13}\}),
((0.5/(x_7,x_1),\{u_{11}\}),(0.3/(x_8,x_1),\{u_8,\\&u_{12}\}),
(0.7(x_7,x_2),\{u_{14}\}), (0.3/(x_8,x_2),\{u_8,u_{14}\}),
(0.3/(x_3,x_3),\{u_1,\\&u_5, u_6,u_9,u_{10}, u_{13}\}),
 ( 0.3/(x_7,x_3),\{u_2,u_5,u_9,u_{10}\}),
 (0.3/(x_8,x_3),\{u_2,\\&u_5,u_{10},u_{12}\}),
(0.9/(x_3,x_4), \{u_6\}),(0.7/(x_7,x_4),\{u_2,u_6\}),
 (0.3/(x_8,x_4), \\&\{u_2,u_8,u_{12}\}),
(0.6/(x_3,x_5),\{u_6,u_9,u_{13}\}),
(0.6/(x_7,x_5),\{u_2\})\bigg\} \\
\end{array}
$$

\end{exmp}
%------------------------------------------------------------

\begin{prop}
Let $R_{F_1}$ and $R_{F_2}$ be two FP-soft relations. Then
\begin{enumerate}
\item $(R_{F_1}^{-1})^{-1}=R_{F_1}$
\item $R_{F_1}\subseteq R_{F_2}\Rightarrow R_{F_1}^{-1}\subseteq R_{F_2}^{-1}$
\end{enumerate}
\end{prop}
\textbf{Proof:}
\begin{enumerate}
\item $\alpha(R_{F_1}^{-1})^{-1}\beta=\beta R_{F_1}^{-1} \alpha=\alpha R_{F_1}\beta$
\item $\alpha R_{F_1}\beta\subseteq \alpha R_{F_2}\beta\Rightarrow \beta R_{F_1}^{-1}\alpha\subseteq
\beta R_{F_2}^{-1}\alpha\Rightarrow R_{F_1}^{-1}\subseteq R_{F_2}^{-1}$
\end{enumerate}

%---------------------------------------------------------------------------------
\begin{defn}
If $R_{F_1}$ is a fuzzy parametrized soft relation from $\Gamma_X$
to $\Gamma_Y$ and $R_{F_2}$ is a fuzzy parametrized soft relation
from $\Gamma_Y$ to $\Gamma_Z$, then a composition of two FP-soft
relations $R_{F_1}$ and $R_{F_2}$ is defined by
$$\alpha (R_{F_1}\circ R_{F_2})\gamma =(\alpha R_{F_1}\beta) \wedge (\beta R_{F_2}\gamma) $$
\end{defn}
\begin{prop}
Let $R_{F_1}$ and $R_{F_2}$ be two FP-soft relation from $\Gamma_X$
to $\Gamma_Y$. Then, $(R_{F_1}\circ R_{F_2})^{-1}=R_{F_2}^{-1}\circ
R_{F_1}^{-1}$
\end{prop}
\textbf{Proof:}
$$
\begin{array}{rl}
\alpha (R_{F_1}\circ R_{F_2})^{-1})\gamma&=\gamma (R_{F_1}\circ
R_{F_2})\alpha\\&=(\gamma R_{F_1}\beta) \wedge(\beta
R_{F_2}\alpha)\\&= (\beta R_{F_2}\alpha)\wedge (\gamma
R_{F_1}\beta)\\&=(\alpha R_{F_2}^{-1}\beta)\wedge (\beta
R_{F_1}^{-1}\gamma)\\&=\alpha(R_{F_2}^{-1}\circ R_{F_1}^{-1})\gamma \\
\end{array}
$$
%-------------
Therefore we obtain

$(R_{F_1}\circ R_{F_2})^{-1}=R_{F_2}^{-1}\circ R_{F_1}^{-1}$

%------------------------------------------------------------
\begin{defn}
An FP-soft relation $R_{F}$ on $\Gamma_X$ is said to be an FP-soft
symmetric relation if $\alpha R_{F_1}\beta\Rightarrow \beta
R_{F_1}\alpha, \forall \alpha, \beta \in \Gamma_X$.
\end{defn}
\begin{defn}
An FP-soft relation $R_{F}$ on $\Gamma_X$ is said to be an FP-soft
transitive relation if $ R_{F}\circ R_{F}\subseteq R_{F}$, that is,
$\alpha R_{F}\beta$ and $\beta R_{F}\gamma \Rightarrow \alpha
R_{F}\gamma, \forall \alpha, \beta, \gamma \in \Gamma_X$.
\end{defn}\begin{defn}
An FP-soft relation $R_{F}$ on $\Gamma_X$ is said to be an FP-soft
reflexive relation if $\alpha R_{F}\alpha , \forall \alpha \in
\Gamma_X$.
\end{defn}

%--------------------------------------------------------------
\begin{defn}
An FP-soft relation $R_{F}$ on $\Gamma_X$ is said to be an FP-soft
equivalence relation if it is symmetric, transitive and reflexive.
\end{defn}

%--------------------------------------------------------------

\begin{exmp}\label{7}
Let $U=\{u_1,u_2,u_3,u_4,u_5, u_6,u_7,u_8\}$,
$E=\{x_1,x_2,x_3,x_4,x_5,\\x_6,x_7,x_8\}$ and $X=
\{0.5/x_1,0.7/x_2,0.3/x_3\}$ be a fuzzy subsets over $E$. Suppose
that
$$
\begin{array}{rl}
\Gamma_{X}= & \bigg\{(0.5/x_1,\{u_1, u_3, u_4,u_6,
u_7,u_8\}),(0.7/x_2,\{u_3,
u_7, u_8\}),\\& (0.3/x_3,\{u_1,u_2,u_4,u_5, u_6,u_9\})\bigg\} \\
\end{array}
$$
Then, a cartesian product on $\Gamma_X$ is obtained as follows
$$
\begin{array}{rl}

\Gamma_X\widehat{\times} \Gamma_X=& \bigg\{ (0.5/(x_1,x_1),\{u_1,
u_3, u_4,u_6, u_7,u_8\}) ,\\&(0.5/(x_1,x_2),\{u_3, u_7, u_8\}),
(0.3/(x_1,x_3),\{u_1,u_4, u_6\}),\\&(0.5/(x_2,x_1),\{u_3, u_7,
u_8\}) ,(0.7/(x_2,x_2),\{u_3, u_7,
u_8\}),\\&(0.3/(x_3,x_1),\{u_1,u_4,
u_6\}), (0.3/(x_3,x_3),\{u_1,u_2,u_4,u_5, u_6,u_9\})\bigg\} \\
\end{array}
$$
Then, we get a fuzzy parametrized soft relation $R_F$ on $F_X$ as
follows
$$
\alpha{R_{F}}\beta\Leftrightarrow \mu_{X\widehat{\times}
Y}(x_i,x_j)/(x_i,x_j))\geq0.3\quad(1\leq i,j \leq 3)
$$
Then
$$
\begin{array}{rl}
 {R_F}=& \bigg\{ (0.5/(x_1,x_1),\{u_1, u_3, u_4,u_6, u_7,u_8\})
,(0.5/(x_1,x_2),\{u_3, u_7, u_8\}),\\& (0.3/(x_1,x_3),\{u_1,u_4,
u_6\}),(0.5/(x_2,x_1),\{u_3, u_7, u_8\}),
\\&(0.7/(x_2,x_2),\{u_3, u_7, u_8\}),(0.3/(x_3,x_1),\{u_1,u_4,
u_6\}),\\& (0.3/(x_3,x_3),\{u_1,u_2,u_4,u_5, u_6,u_9\})\bigg\} \\
\end{array}
$$
$R_{F}$ on $\Gamma_X$ is an FP-soft equivalence relation because it
is symmetric, transitive and reflexive.
\end{exmp}
%----------------

\begin{prop}
If $R_{F}$ is symmetric if and only if $R_{F}^{-1}$  is so.
\end{prop}
\textbf{Proof:} If $R_{F}$  is symmetric, then $\alpha
R_{F}^{-1}\beta $ = $\beta R_{F}\alpha$ =$\alpha R_{F}\beta $=$\beta
R_{F}^{-1}\alpha$. So,  $ R_{F}^{-1}$ is symmetric.

Conversely, if $ R_{F}^{-1}$ is symmetric, then $\alpha R_{F}\beta $
= $\alpha (R_{F}^{-1})^{-1}\beta $   = $\beta (R_{F}^{-1})\alpha $ =
$\alpha (R_{F}^{-1})\beta $ =  $\beta R_{F}\alpha $ So, $R_{F}$ is
symmetric.

\begin{prop}$R_{F}$ is symmetric if and only if
$R_{F}^{-1}$=$R_{F}$
\end{prop}
\textbf{Proof:} If $R_{F}$ is symmetric, then $\alpha
R_{F}^{-1}\beta $ = $\beta R_{F}\alpha$ = $ \alpha R_{F}\beta$. So,
$R_{F}^{-1}$ = $R_{F}$.

Conversely, if $R_{F}^{-1}$ =$R_{F}$, then $\alpha R_{F}\beta$ =
$\alpha R_{F}^{-1}\beta$ = $\beta R_{F}\alpha $. So, $R_{F}$ is
symmetric.

\begin{prop}If $R_{F_1}$ and
$R_{F_2}$ are symmetric relations on $\Gamma_X$, then $R_{F_1}\circ
R_{F_2}$ is symmetric on $\Gamma_X$  if and only if $R_{F_1}\circ
R_{F_2}$=$R_{F_2}\circ R_{F_1}$
\end{prop}
\textbf{Proof:} If $R_{F_1}$ and $R_{F_2}$ are symmetric, then it
implies $R_{F_1}^{-1}=R_{F_1}$ and $R_{F_2}^{-1}=R_{F_2}$. We have
$(R_{F_1}\circ R_{F_2})^{-1}$ = $R_{F_2}^{-1}\circ R_{F_1}^{-1}$.
then $R_{F_1}\circ R_{F_2}$ is symmetric. It implies $R_{F_1}\circ
R_{F_2}$ = $(R_{F_1}\circ R_{F_2})^{-1}$ = $R_{F_2}^{-1}\circ
R_{F_1}^{-1}$=$R_{F_2}\circ R_{F_1}$.

Conversely, $(R_{F_1}\circ R_{F_2})^{-1}$ = $R_{F_2}^{-1}\circ
R_{F_1}^{-1}$ = $R_{F_2}\circ R_{F_1}$ = $R_{F_1}\circ R_{F_2}$. So,
$R_{F_1}\circ R_{F_2}$ is symmetric.

\begin{ccorollary}
If $R_{F}$ is symmetric, then $R_{F}^n$  is symmetric for all
positive integer n, where $R_{F}^n=\underbrace{R_{F}\circ
R_{F}\circ...\circ R_{F}}_{n\,\, times}$.
\end{ccorollary}

\begin{prop}If $R_{F}$ is transitive, then $R_{F}^{-1}$ is also transitive.
\end{prop}
\textbf{Proof:}
$$
\begin{array}{rl}\alpha R_{F}^{-1}\beta &= \beta R_{F}\alpha \supseteq
\beta( R_{F}\circ R_{F}) \alpha\\&=(\beta R_{F}\gamma) \wedge
(\gamma R_{F}\alpha) \\&= (\gamma R_{F}\alpha)\wedge (\beta
R_{F}\gamma) \\&= (\alpha R_{F}^{-1}\gamma) \wedge (\gamma
R_{F}^{-1}\beta) \\&=\alpha (R_{F}^{-1} \circ R_{F}^{-1}) \beta
\end{array}
$$
So, $R_{F}^{-1} \circ R_{F}^{-1}\subseteq R_{F}^{-1}$. The proof is
completed.

\begin{prop}
If $R_{F}$ is transitive
then $R_{F}\circ R_{F}$ is so.
\end{prop}
\textbf{Proof:}
$$
\begin{array}{rl}
\alpha(R_{F}\circ R_{F})\beta &= (\alpha R_{F}\gamma) \wedge (\gamma
\frac{\frac{}{}}{}R_{F}\beta)\\&=\alpha (R_{F}\circ R_{F})\gamma
\wedge \gamma(R_{F}\circ R_{F})\beta\\&=\alpha(R_{F}\circ R_{F}\circ
R_{F}\circ R_{F})\beta
\end{array}
$$
So, $\alpha(R_{F}\circ R_{F}\circ R_{F}\circ R_{F})\beta \subseteq
\alpha(R_{F}\circ R_{F})\beta$. The proof is completed.

\begin{prop}
If $R_{F}$ is reflexive then $R_{F}^{-1}$ is so.
\end{prop}
\textbf{Proof:}
 $\alpha R_{F}^{-1}\beta = \beta R_{F}\alpha\subseteq \alpha R_{F}\alpha = \alpha R_{F}^{-1}\alpha$
and $\beta R_{F}^{-1}\alpha=\alpha R_{F}\beta \subseteq \alpha
R_{F}\alpha=\alpha R_{F}^{-1}\alpha$. The proof is completed.

\begin{prop} If $R_{F}$  is symmetric and transitive, then $R_{F}$  is
reflexive.
\end{prop}
\textbf{Proof:} Proof can be made easily by using Definition 4.1,
Definition 4.2 and Definition 4.3.

%------------------------------------------------------------
\begin{defn}
Let $\Gamma_X\in FPS(U)$, $R_F$ be an FP-soft equivalence relation
on $\Gamma_X$ and $\alpha\in R_F$. Then, an equivalence class of
$\alpha$, denoted by $[\alpha]_{R_F}$, is defined as
$$[\alpha]_{R_F}=\{\beta:\alpha R_F\beta\}.$$
\end{defn}
%-----------------------------------------------------------
\begin{exmp}\label{8}Let us consider the Example \ref{7}.
 Then an equivalence class of \\$(x_1,\{u_1, u_3,
u_4,u_6, u_7,u_8\})$ will be as follows.
$$
\begin{array}{rl}
[(0.5/x_1,\{u_1, u_3, u_4,u_6, u_7,u_8\})]_{R_F}=& \bigg\{
(0.5/x_1,\{u_1, u_3, u_4,u_6, u_7,u_8\}),\\&(0.7/x_2,\{u_3, u_7,
u_8\}),(0.3/x_3,\\&\{u_1,u_2,u_4,u_5,
u_6,u_9\})\bigg\} \\
\end{array}
$$
\end{exmp}
%--------------------------------------------------------------
\section{Decision Making Method}
In this section, we construct a soft fuzzification operator and a
decision making method on FP-soft relations.
\begin{defn}
Let $\Gamma_X\in FPS(U)$ and  $ R_F$ be a FP-soft relation on
$\Gamma_X$. Then fuzzification operator, denoted by $s_{R_{F}}$, is
defined by
$$
s_{R_{F}}: R_{F}\to F(U) , \quad s_{R_{F}}(X\times X,U)=
\{\mu_{R_{F}}(u)/u: u\in U\}
$$
where
$$
\mu_{R_{F}} (u)=\frac{1}{|X\times
X|}\sum_{j}\sum_{i}\mu_{R_{F}}(x_i,x_j)\chi(u)
$$
and where
$$
\chi(u)=\left\{
\begin{array}{ll}
1, & u\in f_{R_{F}}(x_i,x_j)\\
0, & u \notin f_{R_{F}}(x_i,x_j)
\end{array}\right.
$$
Note that $|X\times X|$ is the cardinality of $X\times X$.
\end{defn}
%------------------------------------------------------------
Now; we can construct a decision making method on FP-soft relation
by the following algorithm;

\begin{enumerate}
    \item construct a feasible fuzzy subset ${X}$  over $E$,
    \item construct a FP-soft set $\Gamma_X$ over $U$,
    \item construct a FP-soft relation ${R_{F}}$ over $\Gamma_X$ according to the requests,
  \item calculate the fuzzification operator $s_{R_{F}}$ over $R_F$,
 \item select the objects, from $s_{R_{F}}$, which have the largest membership value.
\end{enumerate}

%--------------------------------------------------------------
\begin{exmp}
A customer, Mr. X, comes to the auto gallery agent to buy a car
which is over middle class. Assume that an auto gallery agent has a
set of different types of car $U=\{u_1,u_2,u_3,u_4,u_5,
u_6,u_7,u_8\}$, which may be characterized by a set of parameters
$E=\{x_1,x_2,x_3,x_4\}$. For $i=1,2,3,4$ the parameters $x_i$ stand
for ``safety'', ``cheap'', ``modern'' and ``large'', respectively.
If Mr. X has to consider own set of parameters, then we select a car
on the basis of the set of customer parameters by using the
algorithm as follows.

%-------------------------------------------------------------------
\begin{enumerate}
\item  Mr X  constructs a fuzzy sets $X$ over $E$,

$X= \{0.5/x_1,0.7/x_2,0.3/x_3\}$

\item Mr X  constructs a FP-soft set $\Gamma_X$ over $U$,
$$
\begin{array}{rl}
\Gamma_X=&\{(0.5/x_1,\{u_1, u_3, u_4,u_6, u_7,u_8\}),(0.7/x_2,\{u_3,
u_7, u_8\}), (0.3/x_3,\\&\{u_1,u_2,u_4,u_5, u_6,u_9\}
\end{array}
$$
\item the fuzzy parametrized soft
relation $R_F$ over $\Gamma_X$ is calculated according to the Mr X's
requests (The car must be a over middle class, it means the
membership degrees are over 0.5),
$$
\begin{array}{rl}
R_F =& \bigg\{ (0.5/(x_1,x_1),\{u_1, u_3, u_4,u_6, u_7,u_8\})
,(0.5/(x_1,x_2),\{u_3, u_7,\\& u_8\}),(0.5/(x_2,x_1),\{u_3, u_7,
u_8\}),(0.7/(x_2,x_2),\{u_3, u_7, u_8\})\bigg\}
\end{array}
$$

\item the soft fuzzification operator $s_{R_{F}}$ over $R_F$ is calculated
as follows
$$
\begin{array}{rl}
s_{R_{F}}=& \bigg\{ (0.055/u_1, 0.0/u_2, 0.244/u_3, 0.055/u_4,
0.0/u_5, 0.055/u_6, 0.244/u_7, \\&0.244/u_8\}\bigg\}
\end{array}
$$
\item
now, select the optimum alternative objects $u_3$, $u_7$ and $u_8$
which have the biggest membership degree 0.244 among the others.
\end{enumerate}
\end{exmp}
%--------------------------------------------------------------
\section{Conclusion}
We first gave most of the fundamental definitions of the operations
of fuzzy sets, soft sets and FP-soft sets are presented. We then
defined relations on FP-soft sets and studied some of their
properties. We also defined symmetric, transitive and reflexive
relations on the FP-soft sets. Finally, we construct a decision
making method and gave an application which shows that this method
successfully works. We have used  a t-norm, which is minimum
operator, the above relation. However, application areas the
relations  can be expanded using the above other norms in the
future.

\end{document}